\documentclass[12pt]{report}
\usepackage{amssymb,amsmath,makeidx,verbatim}
\newcommand{\R}{\mathbb{R}}

\newcommand{\N}{\mathbb{N}}
\newcommand{\C}{\mathbb{C}}
\newcommand{\Z}{\mathbb{Z}}

\newcommand{\be}{\begin{enumerate}}
\newcommand{\ee}{\end{enumerate}}
\newcommand{\bq}{\begin{eqnarray*}}
\newcommand{\eq}{\end{eqnarray*}}

\begin{document}
\newcommand{\disp}{\displaystyle}
\thispagestyle{empty}
\begin{center}
\textbf{Series Analysis and Schwartz Algebras of Spherical Convolutions on Semisimple Lie Groups.\\}
\ \\
\textbf{Olufemi O. Oyadare}\\
\ \\
\end{center}
\begin{quote}

\indent {\bf Abstract}. We give the exact contributions of \textit{Harish-Chandra transform,} $(\mathcal{H}f)(\lambda),$ of Schwartz functions $f$ to the harmonic analysis of \textit{spherical convolutions} and the corresponding $L^{p}-$ Schwartz algebras on a connected semisimple Lie group $G$ (with finite center). One of our major results gives the proof of how the \textit{Trombi-Varadarajan Theorem} enters into the spherical convolution transform of $L^{p}-$ Schwartz functions.
\end{quote}
\ \\
\ \\
\ \\
\ \\
\ \\
\textbf{Subject Classification:} $43A85, \;\; 22E30, \;\; 22E46$\\
\textbf{Keywords:} Harish-Chandra Transforms; Semisimple Lie groups; Harish-Chandra's Schwartz algebras\\
\ \\
{\bf $\bf{1}\;\;\;\;$ Introduction}\\
\ \\

Let $G$ be a connected semisimple Lie group with finite center,
and denote the Harish-Chandra-type Schwartz spaces of functions on $G$
by ${\cal C}^p(G)$, $0 < p \leq 2.$ We know that ${\cal C}^p(G)\subset L^p(G)$ for every such $p$, and if $K$ is a maximal
compact subgroup of $G$ such that ${\cal C}^p(G//K)$ represents the subspace
of ${\cal C}^p(G)$ consisting of the $K-$bi-invariant functions, Trombi and Varadarajan ([$9.$]) have shown that the spherical Fourier
transform $f \mapsto \widehat{f}$ is a linear topological isomorphism
of ${\cal C}^p(G//K)$ onto the spaces $\bar{\mathcal{Z}}({\mathfrak{F}}^{\epsilon})$,
$\epsilon = \left(2/p\right)-1,$
\ \\
\ \\
$\overline{\textmd{Department\;of\;Mathematics,}}$
Obafemi Awolowo University,
Ile-Ife, $220005,$ Nigeria.\\
\text{E-mail: \textit{femi\_oya@yahoo.com}}
\ \\
\ \\
consisting of rapidly decreasing functions on certain sets
${\mathfrak{F}}^{\epsilon}$ of elementary spherical functions.\\

We show the existence of a \textit{hyper-function} on both $G$ and ${\mathfrak{F}}^{1}$ (here named a \textit{spherical convolution}) whose restriction to the group identity element, $e,$ coincides with the \textit{spherical Fourier transforms,} $f \mapsto \widehat{f},$ of Schwartz functions $f$ on $G$ and which affords us the opportunity of embarking on a more inclusive harmonic analysis on $G.$ Indeed $[8a.]$ contains a more general Plancherel formula for the collection of these functions. As a function on $G$ its series expansion is in the present paper studied. We show that, aside from the fact that the spherical Fourier transforms, $\widehat{f}(\lambda),$ is the constant term of this series expansion, there is a region in $G$ where the spherical convolution is \textit{essentially} $\widehat{f}(\lambda).$ Various algebras of these functions are thus studied and ultimately embedded in $L^{2}(G).$ It is however clear that the results in $[8.]$ and in the present paper may be extended to include what may be termed as \textit{the Harish-Chandra-type Schwartz spaces of Eisenstein Integrals on $G.$} The author has recently used the idea of a spherical convolution to give an explicit computation of the image of $\mathcal{C}^{2}(G)$ under the Harish-Chandra transform, $[8b.]$ thus giving a concrete realization of the \textit{abstract} results of Arthur, $[2.],$ and showing the direct contribution of the Plancherel formula to Harish-Chandra transform on $G.$\\

The following is the breakdown of each of the remaining sections of the paper. $\S2.$ contains the preliminaries to the research containing the structure theory, spherical functions and Schwartz algebras on $G,$ while the series analysis of spherical convolutions on $G$ is the subject of $\S3.$ The relationship existing among the Schwartz algebras of functions and those of spherical convolutions is considered in $\S4.$\\
\ \\
{\bf $\bf{2}\;\;\;\;$ Preliminaries}\\

For the connected semisimple Lie group
$G$ with finite center, we denote its Lie algebra by $\mathfrak{g}$
whose \textit{Cartan decomposition} is given as $\mathfrak{g} = \mathfrak{t}\oplus\mathfrak{p}.$ Denote by $\theta$ the \textit{Cartan involution} on $\mathfrak{g}$ whose collection of fixed points is $\mathfrak{t}.$
We also denote by $K$ the analytic subgroup of $G$ with Lie
algebra $\mathfrak{t}.$  $K$ is then a maximal compact subgroup of $G.$
Choose a maximal abelian subspace  $\mathfrak{a}$ of $\mathfrak{p}$ with algebraic
dual $\mathfrak{a}^*$ and set $A =\exp \mathfrak{a}.$  For every $\lambda \in \mathfrak{a}^*$ put
$$\mathfrak{g}_{\lambda} = \{X \in \mathfrak{g}: [H, X] =
\lambda(H)X, \forall  H \in \mathfrak{a}\},$$ and call $\lambda$ a restricted
root of $(\mathfrak{g},\mathfrak{a})$ whenever $\mathfrak{g}_{\lambda}\neq\{0\}.$\\

Denote by $\mathfrak{a}'$ the open subset of $\mathfrak{a}$
where all restricted roots are $\neq 0,$ and call its connected
components the \textit{Weyl chambers.}  Let $\mathfrak{a}^+$ be one of the Weyl
chambers, define the restricted root $\lambda$ positive whenever it
is positive on $\mathfrak{a}^+$ and denote by $\triangle^+$ the set of all
restricted positive roots. Members of $\triangle^+$ which form a basis for $\triangle$ and can not be written as a linear combination of other members of $\triangle^+$ are called \textit{simple.} We then have the \textit{Iwasawa
decomposition} $G = KAN$, where $N$ is the analytic subgroup of $G$
corresponding to $\mathfrak{n} = \sum_{\lambda \in \triangle^+} \mathfrak{g}_{\lambda}$,
and the \textit{polar decomposition} $G = K\cdot
cl(A^+)\cdot K,$ with $A^+ = \exp \mathfrak{a}^+,$ and $cl(A^{+})$ denoting the closure of $A^{+}.$\\

If we set $M = \{k
\in K: Ad(k)H = H$, $H\in \mathfrak{a}\}$ and $M' = \{k
\in K : Ad(k)\mathfrak{a} \subset \mathfrak{a}\}$ and call them the
\textit{centralizer} and \textit{normalizer} of $\mathfrak{a}$ in $K,$ respectively, then (see $[5.]$, p. $284$);
(i) $M$ and $M'$ are compact and have the same Lie algebra and
(ii) the factor  $\mathfrak{w} = M'/M$ is a finite group called the \textit{Weyl
group}. $\mathfrak{w}$ acts on $\mathfrak{a}^*_{\C}$ as a group of linear
transformations by the requirement $$(s\lambda)(H) =
\lambda(s^{-1}H),$$ $H \in \mathfrak{a}$, $s \in \mathfrak{w}$, $\lambda \in
\mathfrak{a}^*_\mathbb{\C}$, the complexification of $\mathfrak{a}^*$.  We then have the
\textit{Bruhat decomposition} $$G = \bigsqcup_{s\in \mathfrak{w}} B m_sB$$ where
$B = MAN$ is a closed subgroup of $G$ and $m_s \in M'$ is the
representative of $s$ (i.e., $s = m_sM$). The Weyl group invariant members of a space shall be denoted by the superscript $^{\mathfrak{w}}$ while $\mid\mathfrak{w} \mid$ represents the cardinality of $\mathfrak{w}.$\\

Some of the most important functions on $G$ are the \textit{spherical
functions} which we now discuss as follows.  A non-zero continuous
function $\varphi$ on $G$ shall be called a \textit{(zonal) spherical
function} whenever $\varphi(e)=1,$ $\varphi \in C(G//K):=\{g\in
C(G)$: $g(k_1 x k_2) = g(x)$, $k_1,k_2 \in K$, $x \in G\}$ and $f*\varphi
= (f*\varphi)(e)\cdot \varphi$ for every $f \in C_c(G//K),$ where $(f \ast g)(x):=\int_{G}f(y)g(y^{-1}x)dy.$  This
leads to the existence of a homomorphism $\lambda :
C_c(G//K)\rightarrow \C$ given as $\lambda(f) = (f*\varphi)(e)$.
This definition is equivalent to the satisfaction of the functional relation $$\int_K\varphi(xky)dk = \varphi(x)\varphi(y),\;\;\;x,y\in G.$$\\

It has been shown by Harish-Chandra [$6.$] that spherical functions on $G$
can be parametrized by members of $\mathfrak{a}^*_{\C}.$  Indeed every
spherical function on $G$ is of the form $$\varphi_{\lambda}(x) = \int_Ke^{(i\lambda-p)H(xk)}dk,\; \lambda
\in \mathfrak{a}^*_{\C},$$  $\rho =
\frac{1}{2}\sum_{\lambda\in\triangle^+} m_{\lambda}\cdot\lambda,$ where
$m_{\lambda}=dim (\mathfrak{g}_\lambda),$ and that $\varphi_{\lambda} =
\varphi_{\mu}$ iff $\lambda = s\mu$ for some $s \in \mathfrak{w}.$ Some of
the well-known properties of spherical functions are $\varphi_{-\lambda}(x^{-1}) =
\varphi_{\lambda}(x),$ $\varphi_{-\lambda}(x) =
\bar{\varphi}_{\bar{\lambda}}(x),$ $\mid \varphi_{\lambda}(x) \mid\leq \varphi_{\Re\lambda}(x),$ $\mid \varphi_{\lambda}(x)\mid\leq \varphi_{i\Im\lambda}(x),$ $\varphi_{-i\rho}(x)=1,$ $\lambda \in \mathfrak{a}^*_{\C},$ while $\mid \varphi_{\lambda}(x) \mid\leq \varphi_{0}(x),\;\lambda \in i\mathfrak{a}^{*},\;x \in G.$ Also if $\Omega$ is the \textit{Casimir operator} on $G$ then
$$\Omega\varphi_{\lambda} = -(\langle\lambda,\lambda\rangle +
\langle \rho, \rho\rangle)\varphi_{\lambda},$$ where $\lambda \in
\mathfrak{a}^*_{\C}$ and $\langle\lambda,\mu\rangle
:=tr(adH_{\lambda} \ adH_{\mu})$ for elements $H_{\lambda}$, $H_{\mu}
\in {\mathfrak{a}}.$ This differential equation may be written simply as $\Omega\varphi_{\lambda} = \gamma(\Omega)(\lambda)\varphi_{\lambda},$ where $\lambda \mapsto \gamma(\Omega)(\lambda)$ is the well-known \textit{Harish-Chandra homomorphism.} The elements $H_{\lambda}$, $H_{\mu}
\in {\mathfrak{a}}$  are uniquely defined by the requirement that $\lambda
(H)=tr(adH \ adH_{\lambda})$ and $\mu
(H)=tr(adH \ adH_{\mu})$ for every $H \in {\mathfrak{a}}$ ([$5.$],
Theorem $4.2$). Clearly $\Omega\varphi_0 = 0.$\\

Due to a hint dropped
by Dixmier $[4.]$ $(cf.\;[9.])$ in his discussion of some functional calculus,
it is necessary to recall the notion of
a \textit{`positive-definite'} function and then discuss the situation for
positive-definite spherical functions.  We call a continuous function
$f: G \rightarrow \C$ (algebraically) positive-definite whenever, for all
$ x_1,\dots,x_m $ in $G$ and all $ \alpha_1,\dots,\alpha_m$ in $\C,$ we have $$\sum^m_{i,j=1}\alpha_i\bar{\alpha}_jf(x^{-1}_i x_j) \geq 0.$$  It
can be shown $(cf.\;[5.])$ that $f(e) \geq 0$ and $|f(x)| \leq f(e)$ for every
$x \in G$ implying that the space ${\cal P}$ of all
positive-definite spherical functions on $G$ is a subset of the
space ${\mathfrak{F}}^{1}$ of all bounded spherical functions on $G.$\\

We know, by the Helgason-Johnson theorem ($[7.]$), that $${\mathfrak{F}}^{1}=
\mathfrak{a}^*+iC_{\rho}$$ where $C_{\rho}$ is the convex hull of $\{s\rho: s \in
\mathfrak{w}\}$ in $\mathfrak{a}^*.$ Defining the \textit{involution} $f^*$ of $f$ as $f^*(x) =
\overline{f(x^{-1})}$, it follows that $f = f^*$ for every $f \in
{\cal P}$, and if $\varphi_{\lambda} \in {\cal P}$, then $\lambda$
and $\bar{\lambda}$ are Weyl group conjugate, leading to a realization of $\mathcal{P}$ as a subset of $\mathfrak{w} \setminus \mathfrak{a}^*_{\C}.$  ${\cal P}$ becomes
a locally compact Hausdorff space when endowed with the \textit{weak $^{*}-$topology} as a subset of $L^{\infty}(G)$.\\

Let $$\varphi_0(x):= \int_{K}\exp(-\rho(H(xk)))dk$$ be denoted
as $\Xi(x)$ and define $\sigma: G \rightarrow \C$ as
$$\sigma(x) = \|X\|$$ for every $x = k\exp X \in G,\;\; k \in K,\; X
\in \mathfrak{a},$ where $\|\cdot\|$ is a norm on the finite-dimensional
space $\mathfrak{a}.$ These two functions are spherical functions on
$G$ and there exist numbers $c,d$ such that $$1 \leq \Xi(a)
e^{\rho(\log a)} \leq c(1+\sigma(a))^d.$$ Also there exists $r
> 0$ such that $c =: \int_G\Xi(x)^2(1+\sigma(x))^{r}dx
< \infty$ ($[11.],$ p. $231$). For each
$0 \leq p \leq 2$ define ${\cal C}^p(G)$ to be the set consisting of
functions $f$ in $C^{\infty}(G)$ for which $$\mu_{a,
b;r}(f) :=\sup_G[|f(a; x ; b)|\Xi (x)^{-2/p}(1+\sigma(x))^r] <
\infty$$ where $a,b \in \mathfrak{U}(\mathfrak{g}_{\C}),$ the \textit{universal
enveloping algebra} of $\mathfrak{g}_{\C},$ $r \in \Z^+, x \in G,$
$f(x;b) := \left.\frac{d}{dt}\right|_{t=0}f(x\cdot(\exp tb))$
and $f(a;x) :=\left.\frac{d}{dt}\right|_{t=0}f((\exp
ta)\cdot x).$ We call ${\cal C}^p(G)$ the Schwartz space on $G$
for each $0 < p \leq 2$ and note that ${\cal C}^2(G)$ is the
well-known (see $[1.]$) Harish-Chandra space of \textit{rapidly decreasing functions} on
$G.$ The inclusions $$C^{\infty}_{c}(G) \subset {\cal C}^p(G)
\subset L^p(G)$$ hold and with dense images. It also follows that
${\cal C}^p(G) \subseteq {\cal C}^q(G)$ whenever $0 \leq p \leq q
\leq 2.$ Each ${\cal C}^p(G)$ is closed under \textit{involution} and the
\textit{convolution}, $*.$ Indeed ${\cal C}^p(G)$ is a Fr$\acute{e}$chet algebra ($[10.],$ p. $69$). We endow ${\cal C}^p(G//K)$
with the relative topology as a subset of ${\cal C}^p(G).$\\

We shall say a function $f$ on $G$ satisfies a \textit{general strong inequality} if for any $r\geq0$ there is a constant $c=c_{r}>0$ such that
$$\mid f(y) \mid \leq c_{r} \Xi(y^{-1}x) (1+\sigma(y^{-1}x))^{-r}\;\;\;\;\;\forall\;x,y \in G.$$ We observe that if $x=e$ then, using the fact that $\Xi(y^{-1})=\Xi(y)$ and $\sigma(y^{-1})=\sigma(y),\;\forall\;y \in G,$ such a function satisfies $$\mid f(y) \mid \leq c_{r} \Xi(y^{-1}) (1+\sigma(y^{-1}))^{-r}=c_{r} \Xi(y) (1+\sigma(y))^{-r},\;\forall\;y \in G,$$ showing that a function on $G$ which satisfies a general strong inequality satisfies in particular a \textit{strong inequality} (in the classical sense of Harish-Chandra, $[11.]$). Members of $\mathcal{C}^{2}(G)=:\mathcal{C}(G)$ are those functions $f$ on $G$ for which $f(g_1; \cdot ; g_2)$ satisfies the strong inequality, for all $g_1,g_2 \in \mathfrak{U}(\mathfrak{g}_{\C}).$ We may then define $\mathcal{C}^{(x)}(G)$ to be those functions $f$ on $G$ for which $f(g_1; \cdot ; g_2)$ satisfies the general strong inequality, for all $g_1,g_2 \in \mathfrak{U}(\mathfrak{g}_{\C})$ and a fixed $x \in G.$ It is clear that $\mathcal{C}^{(e)}(G)=\mathcal{C}(G)$ and that $\bigcup_{x \in G}\mathcal{C}^{(x)}(G),$ which contains $\mathcal{C}(G),$ may be given an inductive limit topology. The seminorms defining this topology will be explicitly given in $\S 4.$\\

For any measurable function $f$ on $G$ we define the \textit{spherical Fourier
transform} $\widehat{f}$ as $$\widehat{f}(\lambda) = \int_G f(x)
\varphi_{-\lambda}(x)dx,$$ $\lambda \in \mathfrak{a}^*_{\C}.$ It
is known (see $[3.]$) that for $f,g \in L^1(G)$ we have:\\
\begin{enumerate}
\item [$(i.)$] $(f*g)^{\wedge} = \widehat{f}\cdot\widehat{g}$ on $ {\mathfrak{F}}^{1}$
whenever $f$ (or $g$) is right - (or left-) $K$-invariant; \item
[$(ii.)$] $(f^*)^{\wedge}(\varphi) =
\overline{\widehat{f}(\varphi^*)}, \varphi \in {\mathfrak{F}}^{1}$; hence
$(f^*)^{\wedge} = \overline{\widehat{f}}$ on ${\cal P}:$ and, if we
define $f^{\#}(g) := \int_{K\times K}f(k_1xk_2)dk_1dk_2,  x\in
G,$ then \item [$(iii.)$] $(f^{\#})^{\wedge}=\widehat{f}$ on ${\mathfrak{F}}^{1}.$
\end{enumerate}

We shall denote the \textit{spherical Fourier transform} $\widehat{f}(\lambda)$ of $f\in \mathcal{C}(G)$ by $(\mathcal{H}f)(\lambda)$ and refer to it as the \textit{Harish-Chandra transforms} of $f.$ Its major properties are well-known and may be found in $[9.].$ It should be noted that $(\mathcal{H}f)(\lambda)=\widehat{f}(\lambda)=\int_G f(y) \varphi_{-\lambda}(y)dy=\int_G f(y) \varphi_{\lambda}(y^{-1})dy=\int_G f(y) \varphi_{\lambda}(y^{-1}e)dy$\\ $=(f\ast\varphi_{\lambda})(e).$ That is, the Harish-Chandra transforms of $f$ is the restriction of the function $$x \mapsto (f\ast\varphi_{\lambda})(x)=:s_{\lambda,f}(x)$$ on $G$ to the identity element. It is therefore worthwhile to explore $s_{\lambda,f}(x)$ in some details for all $x \in G$ in order to put its behaviour at $x=e$ (as the Harish-Chandra transforms of $f$) in a proper and larger perspective.\\

The beauty of studying the entirety of the function $s_{\lambda,f}(x),$ for $\lambda \in \mathfrak{a}^{*}_{\C},\;f \in \mathcal{C}^p(G),\;x \in G,$ which we shall explore in this paper, is that it could be viewed as a transformation in six ($6$) different ways; As $$(1.)\;\;\;\;\;\;\;\;\;\;x\mapsto k_{1}(\lambda):=s_{\lambda,f}(x),\;\mbox{for any $f \in {\cal C}^p(G)$}$$ and $$(2.)\;\;\;\;\;\;\;\;\;\;\;\;\;\;\;x\mapsto k_{2}(f):=s_{\lambda,f}(x),\;\mbox{for any $\lambda \in \mathfrak{a}^{*}_{\C}$},$$ (from where the Plancherel formula for the space of functions $x\mapsto k_{2}(f)$ has recently been computed in $[8a.]$) both of which are maps on $G;$ or as $$(3.)\;\;\;\;\;\;\;\;\;\;\;\;\;\;\;\;f\mapsto l_{1}(\lambda):=s_{\lambda,f}(x),\;\mbox{for any $x \in G$}$$ (which, at $x=e,$ led Harish-Chandra to the consideration of $f\mapsto(\mathcal{H}f)(\lambda):\;cf.\;[9.]$) and $$(4.)\;\;\;\;\;\;\;\;\;\;\;\;\;\;\;\;f\mapsto l_{2}(x):=s_{\lambda,f}(x),\;\mbox{for any $\lambda \in \mathfrak{a}^{*}_{\C}$},$$ both of which are maps on ${\cal C}^p(G);$ or as $$(5.)\;\;\;\;\;\;\;\;\;\;\;\;\;\;\;\;\lambda\mapsto m_{1}(f):=s_{\lambda,f}(x),\;\mbox{for any $x \in G$}$$ and $$(6.)\;\;\;\;\;\;\;\;\;\;\lambda\mapsto m_{2}(x):=s_{\lambda,f}(x),\;\mbox{for any $f \in {\cal C}^p(G)$},$$ both of which are maps on $\mathfrak{a}^{*}_{\C}.$ Hence the function $x\mapsto s_{\lambda,f}(x)$ may rightly be called an \textit{hyper-function} on $G$ whose major contribution to harmonic analysis would be to \textit{absorb} other known functions of the subject and put their results in \textit{proper perspectives,} as we shall establish here for the \textit{Harish-Chandra transform.}\\

In order to know the image of the spherical Fourier transform when
restricted to ${\cal C}^p(G//K)$ we need the following spaces that are central to the statement
of the well-known result of Trombi and Varadarajan [$9.$]. Let $C_\rho$ be the closed convex hull of the (finite) set $\{s\rho :
s\in \mathfrak{w}\}$ in $\mathfrak{a}^*$, i.e., $$C_\rho =
\left\{\sum^n_{i=1}\lambda_i(s_i\rho) : \lambda_i \geq 0,\;\;\sum^n_{i=1}\lambda_i = 1,\;\;s_i \in \mathfrak{w}\right\}$$ where we recall that, for every
$H \in \mathfrak{a},\;\;(s\rho)(H) = \frac{1}{2} \sum_{\lambda\in\triangle^+}
m_{\lambda}\cdot\lambda (s^{-1}H).$\\

Now for each
$\epsilon > 0$ set ${\mathfrak{F}}^{\epsilon} = \mathfrak{a}^*+i\epsilon
C_\rho.$ Each ${\mathfrak{F}}^{\epsilon}$ is convex in $\mathfrak{a}^*_{\C}$ and
$$int({\mathfrak{F}}^{\epsilon}) =
\bigcup_{0<\epsilon'<\epsilon}{\mathfrak{F}}^{\epsilon^{'}}$$
([$9.$], Lemma $(3.2.2)$).  Let us define $\mathcal{Z}({\mathfrak{F}}^{0}) = \mathcal{S}
(\mathfrak{a}^*)$ and, for each $\epsilon>0,$ let
$\mathcal{Z}({\mathfrak{F}}^{\epsilon})$ be the space of all $\C$-valued
functions $\Phi$ such that  $(i.)$ $\Phi$ is defined and holomorphic
on $int({\mathfrak{F}}^{\epsilon}),$ and $(ii.)$ for each holomorphic
differential operator $D$ with polynomial coefficients we have $\sup_{int({\mathfrak{F}}^{\epsilon})}|D\Phi| < \infty.$\\

The space $\mathcal{Z}({\mathfrak{F}}^{\epsilon})$ is converted to a Fr$\acute{e}$chet algebra by equipping it with the
topology generated by the collection, $\| \cdot \|_{\mathcal{Z}({\mathfrak{F}}^{\epsilon})},$ of seminorms given by $\|\Phi\|_{\mathcal{Z}({\mathfrak{F}}^{\epsilon})} := \sup_{int({\mathfrak{F}}^{\epsilon})}|D\Phi|.$  It is known that $D\Phi$ above extends to
a continuous function on all of ${\mathfrak{F}}^{\epsilon}$
([$9.$], pp. $278-279$). An appropriate subalgebra of
$\mathcal{Z}({\mathfrak{F}}^{\epsilon})$ for our purpose is the closed
subalgebra $\bar{\mathcal{Z}}({\mathfrak{F}}^{\epsilon})$ consisting of
$\mathfrak{w}$-invariant elements of $\mathcal{Z}({\mathfrak{F}}^{\epsilon})$,
$\epsilon \geq 0.$ The following (known as the \textit{Trombi-Varadarajan Theorem}) is the major result of $[9.]:$ \textit{Let $0 < p \leq 2$ and
set $\epsilon = \left(2/p\right)-1.$ Then the
spherical Fourier transform $f \mapsto \widehat{f}$ is a linear
topological algebra isomorphism of ${\cal C}^p(G//K)$ onto $\bar{\mathcal{Z}}
({\mathfrak{F}}^{\epsilon}).$} That is, the topological algebra $\bar{\mathcal{Z}}({\mathfrak{F}}^{\epsilon})$ is an isomorphic copy or a realization of ${\cal C}^p(G//K).$\\

In order to find other isomorphic copies or realizations of ${\cal C}^{p}(G//K)$ under the more inclusive general transformation map $$f\mapsto l_{1}(\lambda):=s_{\lambda,f}(x),\;\mbox{for any $x \in G$},$$ we shall now introduce a more general algebra, $\bar{\mathcal{Z}}_{G}({\mathfrak{F}}^{\epsilon}),$ of $\C-$valued functions on $int({\mathfrak{F}}^{\epsilon})\times G$ which, when restricted to $int({\mathfrak{F}}^{\epsilon})\times \exp(N_{0}),$ coincides with $\bar{\mathcal{Z}}({\mathfrak{F}}^{\epsilon}).$ The form of this new algebra is suggested by Theorem $3.5.$ Set $\mathcal{Z}_{G}({\mathfrak{F}}^{0})=\mathcal{S}(\mathfrak{a}^*)\times G$ and let $\mathcal{Z}_{G}({\mathfrak{F}}^{\epsilon}),\;\epsilon>0,$ be the collection of all $\C-$valued functions $\Psi$ ($(\lambda,x)\mapsto \Psi(\lambda,x),\;\forall\;(\lambda,x) \in int({\mathfrak{F}}^{\epsilon})\times G$) such that\\

$(i.)$ $\Psi$ is holomorphic in the variable $\lambda,$ analytic in $x$ and spherical on $G;$\\

$(ii.)$ $\sup_{int({\mathfrak{F}}^{\epsilon})}|D_{1}\Psi| < \infty$ and $\sup_{G}|\Psi D_{2}| < \infty,$ for every holomorphic differential operator $D_{1}$ with polynomial coefficients and every left-invariant differential operator $D_{2}$ on $G$ and\\

$(iii.)$ the restriction of $\Psi$ to $int({\mathfrak{F}}^{\epsilon})\times \{e\}$ (or to $int({\mathfrak{F}}^{\epsilon})\times \exp(N_{0}(A^{+})),$ for some zero neighbourhood $N_{0}(A^{+})$ in $\mathfrak{g},$ as will later be seen in Theorem $3.5$) is (a non-zero constant multiple of) the Harish-Chandra transform, $(\mathcal{H}f)(\lambda)=\hat{f}.$\\

It may be shown, in exact manner as for $\mathcal{Z}({\mathfrak{F}}^{\epsilon})$ above, that the space $\mathcal{Z}_{G}({\mathfrak{F}}^{\epsilon})$ is converted to a Fr$\acute{e}$chet algebra by equipping it with the
topology generated by the collection, $\| \cdot \|_{\mathcal{Z}_{G}({\mathfrak{F}}^{\epsilon})},$ of seminorms given by $$\|\Psi\|_{\mathcal{Z}_{G}({\mathfrak{F}}^{\epsilon})} := \sup_{int({\mathfrak{F}}^{\epsilon})\times G}|D_{1}\Psi D_{2}|.$$ An appropriate subalgebra of $\mathcal{Z}_{G}({\mathfrak{F}}^{\epsilon})$ for our purpose is the closed
subalgebra $\bar{\mathcal{Z}}_{G}({\mathfrak{F}}^{\epsilon})$ consisting of
$\mathfrak{w}$-invariant elements of $\mathcal{Z}_{G}({\mathfrak{F}}^{\epsilon})$,
$\epsilon \geq 0.$ By the time Theorem $3.5$ is established it will be clear that $\bar{\mathcal{Z}}_{\{x\}}({\mathfrak{F}}^{\epsilon})\simeq\bar{\mathcal{Z}}({\mathfrak{F}}^{\epsilon}),$ for every $x$ in some zero neighbourhood $N_{0}(A^{+})$ in $\mathfrak{g}.$ In particular, $\bar{\mathcal{Z}}_{\{e\}}({\mathfrak{F}}^{\epsilon})\simeq\bar{\mathcal{Z}}({\mathfrak{F}}^{\epsilon}).$
\ \\
\ \\
{\bf $\bf{3}\;\;\;\;$ Series Analysis of Spherical Convolutions}\\

Let $f \in \mathcal{C}(G)$ and $\lambda \in \mathfrak{a}^{\ast}_{\C},$ we recall from $[8a.]$ the definition of \textit{spherical convolutions,} $s_{\lambda,f},$ on $G$ corresponding to the pair $(\lambda,f)$ as $$s_{\lambda,f}(x):=(f\ast\varphi_{\lambda})(x),\;\;x \in G.$$ We already know that $s_{\lambda,f}(e)=(\mathcal{H}f)(\lambda),$ where $e$ is the identity element of $G$ and $\lambda \in i\mathfrak{a}^{\ast}.$ This relation between a function on $G$ at the identity element and another function on $i\mathfrak{a}^{\ast}$ suggests we study the full contribution of the Harish-Chandra transforms, $(\mathcal{H}f)(\lambda),$ of $f$ to the properties of $x \mapsto s_{\lambda,f}(x)$ and to seek other functions on $i\mathfrak{a}^{\ast}$ which have not been known in the harmonic analysis of $G,$ but still contribute to a deeper understanding of the structure of $G.$\\

In order to explore the nature of this idea we consider opening up the spherical convolutions $x \mapsto s_{\lambda,f}(x)$ via its \textit{Taylor's series expansion.}\\

{\bf Lemma 3.1.} \textit{Let $N_{0}$ be a neighbourhood of origin in $\mathfrak{g}$ and $t$ be sufficiently small in $\R$ (say $0\leq t\leq 1$). Then $$s_{\lambda,f}(x\exp tX)=\sum^{\infty}_{n=0}\frac{t^{n}}{n!}[\tilde{X}^{n}s_{\lambda,f}](x),$$  where for every $\;X \in N_{0}$ we set $[\tilde{X}^{n}s_{\lambda,f}](x)=\frac{d^{n}}{du^{n}}s_{\lambda,f}(x\exp uX)_{|_{u=0}}$}

{\bf Proof.} The proof follows from a direct application of \textit{Taylor's series expansion,} $[5.],\;p.\;105.\;\;\Box$\\

At $x=e$ and $t=1$ the formula in the Lemma becomes $$s_{\lambda,f}(\exp X)=
\sum^{\infty}_{n=0}\frac{1}{n!}[\tilde{X}^{n}s_{\lambda,f}](e)=s_{\lambda,f}(e)+\sum^{\infty}_{n=1}\frac{1}{n!}[\tilde{X}^{n}s_{\lambda,f}](e)$$ $$=(\mathcal{H}f)(\lambda)+\sum^{\infty}_{n=1}\frac{1}{n!}[\tilde{X}^{n}s_{\lambda,f}](e),\;\;X \in N_{0}.$$ This observation leads quickly to the following result which gives the exact contribution of the Harish-Chandra transforms to the study of spherical convolutions.\\

{\bf Lemma 3.2.} \textit{The Harish-Chandra transforms, $\lambda\mapsto (\mathcal{H}f)(\lambda),\;f \in \mathcal{C}(G),$ is the constant term in the (Taylor's) series expansion of spherical convolutions, $x\mapsto s_{\lambda,f}(x)$ around $x=e,$ for every $\lambda \in \mathfrak{a}^{*}.\;\;\Box$}\\

It may be deduced, from the expansion leading to the proof Lemma $3.2,$ that the only time the remaining terms in $s_{\lambda,f}(\exp X),$ after the (\textit{non-zero}) constant term $(\mathcal{H}f)(\lambda),$ could vanish is when the differential operator $\tilde{X}=0.$ That is, when $X=0.$ It therefore follows that the well-known (Harish-Chandra) harmonic analysis on $G$ $([1.],[2.],[9.]\;\mbox{and}\;[11.])$ has always been that of the consideration of the map $X\mapsto s_{\lambda,f}(\exp X)$ at only $X=0,$ which is the origin of $\mathfrak{g}$ or which corresponds to the identity point of $\exp(\mathfrak{g}).$ Hence, since the constant term, $(\mathcal{H}f)(\lambda),$ of $s_{\lambda,f}(\exp X)$ corresponds indeed to the consideration of the constant term in the asymptotic expansion of (zonal) spherical functions, $\varphi_{\lambda},$ it also follows that other terms in the expansion of $\varphi_{\lambda}$ may be needed to completely understand $f\mapsto s_{\lambda,f}(x).$\\

The expression for $s_{\lambda,f}(\exp X)$ therefore suggests that a \textit{full} harmonic analysis of $G$ may be attained from a close study of the remaining contributions of the \textit{transform} of $f$ given as $$\lambda \longmapsto \frac{t^{n}}{n!}[\tilde{X}^{n}s_{\lambda,f}](x),$$ for all $X \in N_{0},\;n \in \N\cup\{0\},\;x \in G,\; f \in \mathcal{C}(G)$ and sufficiently small values of $t,$ in the same manner that its constant term, $$\lambda \longmapsto (\mathcal{H}f)(\lambda)$$ had been considered.\\

However before considering the transformational properties of spherical convolutions we note the following lemmas which lead to a more inclusive view of the Trombi-Varadarajan Theorem and prepares the ground for its generalization.\\

{\bf Lemma 3.3.} \textit{Let $N_{0}$ be a neighbourhood of origin in $\mathfrak{g},$ $\lambda \in \mathfrak{a}^{\ast}_{\C}$ and $t$ be sufficiently small in $\R$ (say $0\leq t\leq 1$). Then $$s_{\lambda,f}(x\exp tX)=[\sum^{\infty}_{n=0}\frac{t^{n}}{n!}\gamma(\frac{d^{n}}{du^{n}})(\lambda)_{|_{u=0}}]\cdot s_{\lambda,f}(x),$$  for every $\;X \in N_{0},\;x \in G,\;f \in \mathcal{C}(G).$}

{\bf Proof.} We note here that $$[\tilde{X}s_{\lambda,f}](x)=\frac{d}{du}s_{\lambda,f}(x \exp uX)_{|_{u=0}}=\frac{d}{du}(f\ast\varphi_{\lambda})(x \exp uX)_{|_{u=0}}$$ $$=(f\ast\frac{d}{du}\varphi_{\lambda})(x \exp uX)_{|_{u=0}}=\gamma(\frac{d}{du})(\lambda)\cdot(f\ast\varphi_{\lambda})(x \exp uX)_{|_{u=0}}.$$ Hence $$[\tilde{X}^{n}s_{\lambda,f}](x)=\gamma(\frac{d^{n}}{du^{n}})(\lambda)_{|_{u=0}}\cdot(f\ast\varphi_{\lambda})(x \exp uX)_{|_{u=0}}=\gamma(\frac{d^{n}}{du^{n}})(\lambda)_{|_{u=0}}\cdot s_{\lambda,f}(x).\;\;\Box$$

The particular case of setting $x=e$ and $t=1$ in Lemma $3.3$ introduces the Harish-Chandra transforms, $(\mathcal{H}f)(\lambda),$ into the analysis of this series, proving the following.\\

{\bf Lemma 3.4.} \textit{Let $N_{0}$ be a neighbourhood of origin in $\mathfrak{g},$ $f \in \mathcal{C}(G)$ and $\lambda \in \mathfrak{a}^{\ast}.$ Then the spherical convolution function, $x \mapsto s_{\lambda,f}(x)$ is a non-zero constant multiple of the Harish-Chandra transforms, $(\mathcal{H}f)(\lambda),$ on $\exp (N_{0}).$}

{\bf Proof.} Set $x=e$ and $t=1$ into Lemma $3.3$ to have $$s_{\lambda,f}(\exp X)=[\sum^{\infty}_{n=0}\frac{1}{n!}\gamma(\frac{d^{n}}{du^{n}})(\lambda)_{|_{u=0}}]\cdot s_{\lambda,f}(e)=[\sum^{\infty}_{n=0}\frac{1}{n!}\gamma(\frac{d^{n}}{du^{n}})(\lambda)_{|_{u=0}}]\cdot (\mathcal{H}f)(\lambda),$$ with $\sum^{\infty}_{n=0}\frac{1}{n!}\gamma(\frac{d^{n}}{du^{n}})(\lambda)_{|_{u=0}}
=1+[\sum^{\infty}_{n=1}\frac{1}{n!}\gamma(\frac{d^{n}}{du^{n}})(\lambda)_{|_{u=0}}]\neq0.\;\;\Box$\\

Let us denote the non-zero constant in Lemma $3.4$ above by $\kappa.$ The following theorem is a consequence of normalizing the spherical convolutions in Lemma $3.4.$\\

{\bf Theorem 3.5.} (\textbf{Trombi-Varadarajan Theorem for Spherical Convolutions}) \textit{Let $0 < p \leq 2,$ set $\epsilon = \left(2/p\right)-1$ and $x \in \exp(N_{0}).$ Set $\widehat{f}_{x}(\lambda)=\frac{1}{\kappa}s_{\lambda,f}(x)$ for $f \in {\cal C}^p(G//K).$ Then the spherical convolution transforms $f \mapsto \widehat{f}_{x}$ is a linear topological algebra isomorphism of ${\cal C}^p(G//K)$ onto $\bar{\mathcal{Z}} ({\mathfrak{F}}^{\epsilon}).\;\;\Box$}\\

We recover the Trombi-Varadarajan Theorem for Harish-Chandra transforms by setting $x=e$ in Theorem $3.5.$ Indeed, Theorem $3.5$ above says that every $x \in \exp(N_{0})$ (and not just $x=e$) gives a topological algebra isomorphism between ${\cal C}^p(G//K)$ and $\bar{\mathcal{Z}}({\mathfrak{F}}^{\epsilon}).$ However if $x \in G\setminus\exp(N_{0}),$ for any neighborhood $N_{0}$ of zero in $\mathfrak{g},$ Trombi-Varadarajan Theorem may not be appropriate and it may be necessary to seek a more general realization of ${\cal C}^p(G//K)$ under the map $f\mapsto l_{1}(\lambda):=s_{\lambda,f}(x),\;\mbox{for any $x \in G$}.$ Before considering another major result of this paper, giving the fine structure of spherical convolution functions, we state a result on the finiteness of a central integral usually used in the estimation of many other integrals of harmonic analysis on semisimple Lie groups.\\

To this end we define, for every $x \in G,$ the function $x\mapsto d(x)$ as $$d(x)=\int_{G}\Xi^{2}(y^{-1}x)(1+\sigma(y^{-1}x))^{-r}dy.$$ We observe here that $$d(e)=\int_{G}\Xi^{2}(y^{-1})(1+\sigma(y^{-1}))^{-r}dy=\int_{G}\Xi^{2}(y)(1+\sigma(y))^{-r}dy,$$ which is a constant whose proof of finiteness may be found in $[11.],\;p.\;231.$ This constant is crucial to all harmonic analysis of $\mathcal{C}(G)$ and, in particular, to the embedding of $\mathcal{C}(G)$ in $L^{2}(G).$ It is therefore important to understand the nature of $d(x)$ for all $x \in G$ in order to employ it in a more inclusive harmonic analysis on $G.$ We consider the nature of this integral in the following.\\

{\bf Lemma 3.6.} \textit{Let $x \in G.$ Then there exist $r\geq0$ such that $$d(x)=\int_{G}\Xi^{2}(y^{-1}x)(1+\sigma(y^{-1}x))^{-r}dy<\infty.$$}

{\bf Proof.} We already know that $\Xi(y^{-1}x)\leq1.$ Also $$1+\sigma(y^{-1}x)\leq(1+\sigma(y^{-1}))(1+\sigma(x))=(1+\sigma(y))(1+\sigma(x)).$$ It follows therefore that $$d(x)\leq\int_{G}(1+\sigma(y^{-1}x))^{-r}dy\leq(1+\sigma(x))\int_{G}(1+\sigma(y))dy.$$ The last integral in the above inequality is finite if we embark on its computation via the polar decomposition, $G=K\cdot cl(A^{+})\cdot K,$ of $G.\;\;\Box$\\

{\bf Theorem 3.7.} \textit{Let $N_{0}$ be a neighbourhood of origin in $\mathfrak{g}$ where $f$ is a measurable function on $G$ which satisfies the general strong inequality. The integral defining the spherical convolution function, $x \mapsto s_{\lambda,f}(x),$ is absolutely and uniformly convergent for all $x \in \exp(N_{0}),\;\lambda \in i \mathfrak{a}^{\ast}.$ Moreover the transforms $\lambda \mapsto s_{\lambda,f}(x)$ of $f,$ with $\;x \in \exp(N_{0}),$ is a continuous function on $i\mathfrak{a}^{\ast}.$ If $r\geq0$ is such that $d(x)=\int_{G}\Xi^{2}(y^{-1}x)(1+\sigma(y^{-1}x))^{-r}dy<\infty,\;x \in G,$ then $$\mid s_{\lambda,f}(x)\mid\leq d(x)\cdot\mu_{1,1,r}(f),\;\;x \in G,\;\lambda \in i \mathfrak{a}^{\ast}.$$}

{\bf Proof.} We recall that $\mid\varphi_{\lambda}(x) \mid\leq \varphi_{0}(x)=\Xi(x),\;x \in G,\;\lambda \in i\mathfrak{a}^{\ast}.$ Hence $$\mid (f \ast\varphi_{\lambda})(x) \mid \leq \int_{G}\mid f(y)\varphi_{\lambda}(y^{-1}x)\mid dy\leq \mu_{1,1,r}(f)\int_{G}\Xi^{2}(y^{-1}x)(1+\sigma(y^{-1}x))^{-r}dy$$ $=d(x)\cdot\mu_{1,1,r}(f).$ Continuity follows from the use of the Lebesgue's dominated convergence theorem$.\;\;\Box$\\

The following well-known result on the foundational properties of the Harish-Chandra transforms, $\lambda \mapsto (\mathcal{H}f)(\lambda),\;\lambda \in i\mathfrak{a}^{\ast},$ now follows from the general outlook given by Theorem $3.7.$\\

{\bf Corollary 3.8.} $([9.])$ \textit{Let $f$ be a measurable function on $G$ which satisfies the strong inequality. The integral defining the Harish-Chandra transforms, $$(\mathcal{H}f)(\lambda)=\int_{G}f(x)\varphi_{\lambda}(x)dx,$$ is absolutely and uniformly convergent for all $\lambda \in i\mathfrak{a}^{\ast}$ and is continuous on $i\mathfrak{a}^{\ast}.$ If $r\geq0$ is such that $d=\int_{G}\Xi^{2}(y)(1+\sigma(y))^{-r}dy<\infty,$ then $$(\mathcal{H}f)(\lambda)\mid\leq d\mu_{1,1,r}(f),\;\;\lambda \in i \mathfrak{a}^{\ast}.$$}

{\bf Proof.} Set $X=0$ in Theorem $3.7$ to have the first results. The inequality follows if we set $x=e$ and observe that $d(e)=\int_{G}\Xi^{2}(y^{-1})(1+\sigma(y^{-1}))^{-r}dy=d.\;\;\Box$
\ \\

We now consider the image of ${\cal C}^p(G//K)$ under the \textit{full} spherical convolution map, $f\mapsto l_{1}(\lambda):=s_{\lambda,f}(x),\;\mbox{for any $x \in G$}.$ In order to discuss this we have two options. One of the options is to introduce wave-packet that will still have its domain as $\bar{\mathcal{Z}}({\mathfrak{F}}^{\epsilon})$ while using an appropriate Plancherel measure on $\mathfrak{F^{\epsilon}}.$ This option has been explored in $[8a.],\;p.\;34,$ where the $L^{2}$ Plancherel measure, $d\zeta_{x,\lambda}$ on $\mathfrak{F^{1}}$ for the spherical convolution function (when viewed as a function on $G$) was defined to absorb the group variable, $x.$ The results therein suggest that the image of ${\cal C}^p(G//K)$ under the \textit{full} spherical convolution map is indeed possible.\\

The second option is to retain the \textit{spherical Bochner measure,} $d\lambda,$ on (a subset of) ${\mathfrak{F}}^{\epsilon}$ and define the wave-packet as a map on the Fr$\acute{e}$chet algebra $\bar{\mathcal{Z}}_{G}({\mathfrak{F}}^{\epsilon}).$ This will reflect the nature of the full spherical convolution map as a transform of members of $\mathcal{C}^p(G//K)$ whose arguments are (generally) taken from $int({\mathfrak{F}}^{\epsilon})\times G$ (and not just from $int({\mathfrak{F}}^{\epsilon})$ as in the first option).\\

To this end recall the Fr$\acute{e}$chet algebra $\bar{\mathcal{Z}}_{G}({\mathfrak{F}}^{\epsilon}),\;\forall\;\epsilon >0,$ let $\Psi \in \bar{\mathcal{Z}}_{G}({\mathfrak{F}}^{\epsilon})$ and set $$N_{0}(A^{+})=N_{0}\cap A^{+},$$ where $N_{0}$ is a zero neighbourhood in $\mathfrak{g}.$ It is clear that $N_{0}(A^{+})$ is also a zero neighbourhood in $\mathfrak{g}$ and that $\Psi=\Psi(\lambda,x),$ for all $(\lambda,x) \in int({\mathfrak{F}}^{\epsilon})\times G.$ It follows, from Theorem $3.5,$ that $\bar{\mathcal{Z}}_{\{x\}}({\mathfrak{F}}^{\epsilon})\simeq\bar{\mathcal{Z}}({\mathfrak{F}}^{\epsilon}),$ for every $x \in \exp(N_{0}(A^{+})).$ We then have the following.\\

{\bf Lemma 3.9.} \textit{For every $x \in \exp(N_{0}(A^{+}))$ and $\Psi \in \bar{\mathcal{Z}}_{G}({\mathfrak{F}}^{\epsilon}),$ we have that $\Psi(\lambda,x)=\Phi(\lambda),$ for some $\Phi \in \bar{\mathcal{Z}}({\mathfrak{F}}^{\epsilon}).$}\\

We now employ these remarks to define a map from $\bar{\mathcal{Z}}_{G}({\mathfrak{F}}^{\epsilon})$ to ${\cal C}^p(G//K)$ as follows. Let $a \in \bar{\mathcal{Z}}_{G}({\mathfrak{F}}^{\epsilon})$ and $\lambda\mapsto c(\lambda)$ be the Harish-Chandra $c-$function defined on $\mathfrak{F}_{I}:=i\mathfrak{a}^{*}.$ We associate to every $a \in \bar{\mathcal{Z}}_{G}({\mathfrak{F}}^{\epsilon})$ the function $\varphi_{a}$ on $G$ defined as $$\varphi_{a}(x)=\mid\mathfrak{w} \mid^{-1}\int_{\mathfrak{F}_{I}}a(-\lambda,x)\varphi_{-\lambda}(x)c(-\lambda)^{-1}c(\lambda)^{-1}d\lambda,\;\;\;x \in G.$$
It should be noted here that $$\varphi_{a}(x)=\mid\mathfrak{w} \mid^{-1}\int_{\mathfrak{F}_{I}}a(-\lambda,x)\varphi_{-\lambda}(x)c(-\lambda)^{-1}c(\lambda)^{-1}d\lambda$$
$$=\mid\mathfrak{w} \mid^{-1}\int_{\mathfrak{F}_{I}}a(\lambda,x)\varphi_{\lambda}(x)c(\lambda)^{-1}c(-\lambda)^{-1}d(-\lambda)$$
$$=\mid\mathfrak{w} \mid^{-1}\int_{\mathfrak{F}_{I}}a(\lambda,x)\varphi_{\lambda}(x)c(\lambda)^{-1}c(-\lambda)^{-1}d\lambda,$$ which is due to the invarianve of $d\lambda,$ and that $$\varphi_{a}(k_{1}xk_{2})=\varphi_{a}(x),$$ $\forall\;x \in G,\;k_{1},k_{2} \in K,$ being a property inherited from $a$ and $\varphi_{\lambda}.$\\

The (extra) requirement of being spherical on $G$ placed on members of $\bar{\mathcal{Z}}_{G}({\mathfrak{F}}^{\epsilon})$ may at first be seen as a restriction, when compared to the requirements on members of $\bar{\mathcal{Z}}({\mathfrak{F}}^{\epsilon}).$ It however turns out that this extra requirement is what is needed to assure us of the generalization of the \textit{classical} wave-packets (of Trombi-Varadarajan) on $G$ to all of $x \mapsto\varphi_{a}(x).$ This is established as follows.\\

{\bf Lemma 3.10.} \textit{Let $a \in \bar{\mathcal{Z}}_{G}({\mathfrak{F}}^{\epsilon})$ and $N_{0}(A^{+})$ be as defined above. Then, for every $x \in \exp(N_{0}(A^{+})),$ the map $x \mapsto\varphi_{a}(x)$ is the classical wave-packet of $G.$}\\

{\bf Proof.} We observe that, with $\exp tH \in \exp(N_{0}(A^{+})),$ $$a(\lambda,x)=a(\lambda,k_{1}\exp tH k_{2})=a(\lambda,\exp tH)=\Phi(\lambda),$$ for some $\Phi \in \bar{\mathcal{Z}}({\mathfrak{F}}^{\epsilon}).$ Here we have employed the spherical property of $a$ on $G$ in the second equality and Lemma $3.9$ in the third equality$.\;\Box$\\

The above Lemma shows that the definition and properties of the map $x \mapsto\varphi_{a}(x),$ $x \in G,$ is consistent with the relationship (in Lemma $3.4$) existing between spherical convolutions, $s_{\lambda,f}(x)$ and the Harish-Chandra transfroms, $(\mathcal{H}f)(\lambda).$ Hence in order to extend Trombi-Varadarajan Theorem (which gives the image of the algebra $\mathcal{C}^{p}(G//K)$ under $f\mapsto(\mathcal{H}f)(\lambda)$) to all $x \in G$ (under the spherical convolution tranform), it will be necessary to show that $x \mapsto\varphi_{a}(x)$ is the wave-packet of $f\mapsto s_{\lambda,f}(x)$ for all $x\in G.$ According to Lemma $3.10,$ this needs only be done for those $x=k_{1}\exp tHk_{2}$ in $G$ with $\exp tH \notin \exp(N_{0}(A^{+})),$ for any neighbourhood, $N_{0},$ of zero in $\mathfrak{g}.$ We however give a self-contained discussion of these results, the first of which is given below.\\

{\bf Theorem 3.11.} \textit{$\varphi_{a} \in \mathcal{C}^{p}(G//K)$ for every $a \in \bar{\mathcal{Z}}_{G}({\mathfrak{F}}^{\epsilon}).$}\\

In order to finish the establishment of  this Theorem we need some lemmas which give appropriate background for it. Indeed we derive an appropriate bound for $\mid \varphi_{a}(h;u) \mid,$ where $u \in \mathfrak{U}(\mathfrak{g}_{\C})$ and $h$ is well-chosen, and the appropriate collection of seminorms are also in place. These will be considered in a forthcoming paper on \textit{Trombi-Varadarajan Theorem} via the eigenfunction expansion of spherical convolution, which includes the extension of Theorem $3.5$ to all $x \in G.$
\ \\
\ \\
\ \\
{\bf $\bf{4}\;\;\;\;$ Algebras of Spherical Convolutions}\\

We now consider the various algebras of spherical convolutions that have emanated in the course of this research and their relationship with the Harish-Chandra Schwartz algebra, $\mathcal{C}(G),$ on $G$ as well as its distinguished commutative subalgebra, $\mathcal{C}(G//K),$ of (elementary) spherical functions.\\

Define $\mathcal{C}_{\lambda}(G)=\{s_{\lambda,f}:f \in \mathcal{C}(G)\}$ and set $\mathcal{C}_{\lambda,0}(G)=\{s_{\lambda,\varphi_{\lambda}}\},$ for all $\lambda \in \mathfrak{a}^{\ast}_{\C}.$ It is clear that $\bigcup_{\lambda \in \mathfrak{a}^{\ast}_{\C}}\mathcal{C}_{\lambda}(G)$ is contained in $\mathcal{C}(G).$ We may therefore topologize $\bigcup_{\lambda \in \mathfrak{a}^{\ast}_{\C}}\mathcal{C}_{\lambda}(G)$ by giving it the \textit{relative topology} from the topology defined on $\mathcal{C}(G)$ by the seminorms, $\mu_{a,b,r}.$\\

{\bf Lemma 4.1.} \textit{The inclusions $$[\bigcup_{\lambda \in \mathfrak{a}^{\ast}_{\C}}\mathcal{C}_{\lambda,0}(G)]\subset\mathcal{C}(G//K)\subset[\bigcup_{\lambda \in \mathfrak{a}^{\ast}_{\C}}\mathcal{C}_{\lambda}(G)]\subset\mathcal{C}(G)$$ are all proper$.\;\;\Box$}\\

{\bf Theorem 4.2.} \textit{$\bigcup_{\lambda \in \mathfrak{F}^{1}}\mathcal{C}_{\lambda}(G)$ is a closed subalgebra of $\mathcal{C}(G).$}

{\bf Proof.} We recall that $\mu_{a,b;r}(f\ast\varphi_{\lambda})\leq c\mu_{1,b;r+r_{0}}(f)\cdot\mu_{a,1;r}(\varphi_{\lambda}),$ where $c:=\int_{G}\Xi^{2}(x)(1+\sigma(x))^{-r_{0}}dx<\infty$ for some $r_{0}\geq0.$ However $$\mu_{a,1;r}(\varphi_{\lambda})=\sup_{G}[|\varphi_{\lambda}(1; x ; a)|\cdot\Xi (x)^{-1}(1+\sigma(x))^r]$$ $$=\mid \gamma(a)(\lambda)\mid\cdot\sup_{G}[|\varphi_{\lambda}(x)|\cdot\Xi (x)^{-1}(1+\sigma(x))^r]$$
$$\leq M\mid \gamma(a)(\lambda)\mid\cdot\sup_{G}[\Xi (x)^{-1}(1+\sigma(x))^r]<\infty$$ $$(\mbox{since $\varphi_{\lambda}$ is bounded for all $\lambda \in \mathfrak{F}^{1}$}).$$ Hence $\mu_{a,b;r}(f\ast\varphi_{\lambda})<\infty,\;\forall\;\lambda \in \mathfrak{F}^{1}.\;\;\Box$\\

It may be recalled that members of $\mathcal{C}(G)$ are exactly those functions on $G$ whose left and right derivatives satisfy the \textit{strong inequality.} In the light of this observation we define $\mathcal{C}^{(x)}(G)$ as exactly those functions on $G$ whose left and right derivatives satisfy the \textit{general strong inequality,} for each $x \in G.$ Explicitly we set $\mathcal{C}^{(x)}(G)$ as $$\mathcal{C}^{(x)}(G)=\{f:G\mapsto\C:\sup_{y \in G}[|f(a; y ; b)|\cdot\Xi (y^{-1}x)^{-1}(1+\sigma(y^{-1}x))^r] <\infty\},$$ $x \in G.$ A collection of seminorms on each of $\mathcal{C}^{(x)}(G)$ may be given by $$\mu^{(x)}_{a,b;r}(f) :=\sup_{y \in G}[|f(a; y ; b)|\cdot\Xi (y^{-1}x)^{-1}(1+\sigma(y^{-1}x))^r].$$ It is however clear that $\mathcal{C}^{(e)}(G)=\mathcal{C}(G),$ so that $\mathcal{C}(G)\subset\bigcup_{x \in G}\mathcal{C}^{(x)}(G).$\\

{\bf Theorem 4.3.} \textit{The natural inclusion $\bigcup_{x \in G}\mathcal{C}^{(x)}(G)\subset L^{2}(G)$ has a dense image.}

{\bf Proof.} It is known that the natural inclusion of $\mathcal{C}(G)$ in $L^{2}(G)$ has a dense image, $[1.].$ The result therefore follows if we recall that, as sets of functions, $$\mathcal{C}(G)\subset\bigcup_{x \in G}\mathcal{C}^{(x)}(G)\subset L^{2}(G),$$ where the second inclusion holds from the fact that $d(x)<\infty,\;x \in G.\;\;\Box$\\
\ \\

{\bf   References.}
\begin{description}
\item [{[1.]}] Arthur, J.G., \textit{Harmonic analysis of tempered distributions on semisimple Lie groups of real rank one,} Ph.D. Dissertation, Yale University, $1970.$
    \item [{[2.]}] Arthur, J.G., \textit{Harmonic analysis of the Schwartz space of a reductive Lie group,} I. II. (preprint, $1973$).
        \item [{[3.]}] Barker, W.H.,  The spherical Bochner theorem on semisimple Lie groups, \textit{J. Funct. Anal.,} vol. \textbf{20}  ($1975$), pp. $179-207.$
                \item [{[4.]}] Dixmier, J., Op$\acute{\mbox{e}}$rateurs de rang fini dans les repr$\acute{\mbox{e}}$sentations unitaires,\textit{ Publ. math. de l' Inst. Hautes $\acute{\mbox{E}}$tudes Scient.,} tome $\textbf{6}$ ($1960$), pp. $13-25.$
                    \item [{[5.]}] Helgason, S.,  \textit{``Differential Geometry, Lie Groups and Symmetric Spaces,''} Academic Press, New York, $1978.$
                \item [{[6.]}] Helgason, S.,  \textit{``Groups and Geometric Analysis; Integral Geometry, Invariant Differential Operators, and Spherical Functions,''} Academic Press, New York, $1984.$
                \item [{[7.]}] Helgason, S. and Johnson, K.,  The bounded spherical functions on symmetric spaces,  \textit{Advances in Math,} \textbf{3} ($1969$), pp. 586-593.
                        \item [{[8.]}] Oyadare, O. O., ($a.$) On harmonic analysis of spherical convolutions on semisimple Lie groups, \textit{Theoretical Mathematics and Applications,} \textbf{5(3)} ($2015$), pp. $19-36.$ ($b.$) Fourier transform of Schwartz algebras on groups of Harish-Chandra class. Under review.
                    \item [{[9.]}] Trombi, P.C. and Varadarajan, V.S.,  Spherical transforms on semisimple Lie groups, \textit{Ann. of Math.,} \textbf{94} ($1971$), pp. $246-303.$
                \item [{[10.]}] Varadarajan, V.S.,  The theory of characters and the discrete series for semisimple Lie groups, in \textit{Harmonic Analysis on Homogeneous Spaces,} (C.C.  Moore (ed.)) \textit{Proc. of Symposia in Pure Maths.,} vol. \textbf{26} ($1973$), pp. $45-99.$
                    \item [{[11.]}] Varadarajan, V.S., \textit{``An introduction to harmonic analysis on semisimple Lie groups,''}  Cambridge University Press, Cambridge, $1989.$
                        \end{description}

\end{document}